\documentclass[a4paper,12pt]{article}
\advance\textheight3.2cm    
\advance\topmargin-2.0cm
\advance\textwidth2.9cm
\advance\evensidemargin-1.4cm
\advance\oddsidemargin-1.4cm

\begin{document}

\vspace*{-1.5cm}
\centerline{\LARGE \bf OSGA: A fast subgradient algorithm}

\centerline{\LARGE \bf with optimal complexity}

\vspace{1cm}

\centerline{\sl {\large \bf Arnold Neumaier}}

\centerline{\sl Fakult\"at f\"ur Mathematik, Universit\"at Wien}
\centerline{\sl Oskar-Morgenstern-Platz 1, A-1090 Wien, Austria}
\centerline{\sl email: Arnold.Neumaier@univie.ac.at}
\centerline{\sl WWW: http://www.mat.univie.ac.at/\wave neum/}

\vspace{1cm}

\centerline{\today}

\centerline{{\tt osga.tex}}

\vspace{0.5cm}

{\bf Abstract.}
This paper presents an algorithm for approximately minimizing a 
convex function in simple, not necessarily bounded convex domains, 
assuming only that function values and subgradients are available. 
No global information about the objective function is needed apart 
from a strong convexity parameter (which can be put to zero if only 
convexity is known).
The worst case number of iterations needed to achieve a given accuracy 
is independent of the dimension (which may be infinite) and -- apart 
from a constant factor -- best possible under a variety of smoothness 
assumptions on the objective function. 

\begin{flushleft}
\vspace{0.5cm}

{\bf 2010\hspace{.4em} MSC Classification}:
primary 90C25;
secondary 90C60, 49M37, 65K05, 68Q25\\
\end{flushleft}


\section{Introduction}

In the recent years, first order methods for convex optimization have
become prominent again as they are able to solve large-scale problems in
millions of variables (often arising from applications to image 
processing, compressed sensing, or machine learning), where matrix-based
interior point methods cannot even perform a single iteration. 
(However, a matrix-free interior point method by \sca{Fountoulakis} et
al. \cite{FouGZ} works well in some large compressed sensing problems.)

In 1983, \sca{Nemirovsky \& Yudin} \cite{NemY} proved lower bounds on 
the complexity of first order methods (measured in the number of 
subgradient calls needed to achieve a given accuracy) for convex 
optimization under various regularity assumptions for the objective 
functions. (See \sca{Nesterov} \cite[Sections 2.1.2 and 3.2.1]{Nes.intL}
for a simplified account.) They constructed convex, piecewise linear 
functions in dimensions $n>k$, where no first order method can have 
function values more accurate than $O(k^{-1/2})$ after $k$ subgradient 
evaluations. This implies the need for at least $O(\eps^{-2})$ 
subgradient evaluations in the worst case if $f$ is a nondifferentiable 
but Lipschitz continuous  convex function. 
They also constructed convex quadratic functions in dimensions 
$n\ge 2k$ where no first order method can have function values more 
accurate than $O(k^{-2})$ after $k$ gradient evaluations. 
This implies the need for at least $O(\eps^{-1/2})$ gradient evaluations
in the worst case if $f$ is an arbitrarily often differentiable convex 
function. However in case of strongly convex functions with Lipschitz 
continuous gradients, the known lower bounds on the complexity allow
a dimension-independent linear rate of convergence $O(q^k)$ with 
$0<q<1$. 

Algorithms by \sca{Nesterov} \cite{Nes.intL,Nes.smooth,Nes.PD}
(dating back in the unconstrained, not strongly convex case to 1983 
\cite{Nes.dokl}), achieve the optimal complexity order in all three 
cases.
These algorithms need as input the knowledge of global parameters --
a global Lipschitz constant for the objective functions in the 
nonsmooth case, a global Lipschitz constant for the gradient 
in the smooth case, and an explicit constant of strong convexity in 
the strongly convex case. Later many variants were described 
(see, e.g., \sca{Auslender \& Teboulle} \cite{AusT},
\sca{Lan} et al. \cite{LanLM}), some of which are adaptive in the 
sense that they estimate all required constants during the execution 
of the algorithm. \sca{Beck \& Teboulle} \cite{BecT} developed an 
adaptive proximal point algorithm called FISTA, popular in image 
restauration applications. Like all  proximal point based methods, the 
algorithm needs more information about the objective function than
just subgradients, but delivers in return a higher speed of convergence.
\sca{Tseng} \cite{Tse} gives a common uniform derivation of several 
variants of fast first order algorithms based on proximal points.
\sca{Becker} et al.\ \cite{BecCG} (among other thinds) add adaptive 
features to Tseng's class of algorithms, making them virtually 
independent of global (and hence often pessimistic) Lipschitz 
information. Other such adaptive algorithms
include \sca{Gonzaga} et al.\ \cite{GonK,GonKR} and \sca{Meng \& Chen} 
\cite{MenC}. \sca{Devolder} et al.\ \cite{DevGN} show that both the 
nonsmooth case and the smooth case can be understood in a common way 
in terms of inexact gradient methods. 

If the Lipschitz constant is very large, the methods with optimal 
complexity for the smooth case are initially much slower than the
methods that have an optimal complexity for the nonsmooth case.
This counterintuitive situation was remedied by \sca{Lan} \cite{Lan}, 
who provides an algorithm that needs (expensive auxiliary computations 
but) no knowledge about the function except convexity and has the 
optimal complexity, both in the nonsmooth 
case and in the smooth case, without having to know whether or not the
 function is smooth. However, its worst case behavior on 
strongly convex problem is unknown. Similarly, if the constant of 
strong convexity is very tiny, the methods with optimal complexity for 
the strongly convex case are initially much slower than the methods 
that do not rely on strong convexity. Prior to the present work, no 
algorithm was known with optimal complexity both for the general 
nonsmooth case and for the strongly convex case. 

{\bf Content.}
In this paper, we derive an algorithm for approximating a solution 
$\wh x\in C$ of the convex optimization problem
\lbeq{e.f0}
f(\wh x) = \wh f := \min_{x\in C} f(x) 
\eeq
using first order information (function values $f$ and subgradients $g$)
only.
Here $f:C\to\Rz$ is a convex function defined on a nonempty, 
convex subset $C$ of a vector space $V$ with bilinear pairing $\<h,z\>$ 
defined for $z\in V$ and $h$ in the dual space $V^*$. 
The minimum in \gzit{e.f0} exists if there is a point $x^0\in C$ such 
that the level set $\{x\in C\mid f(x)\le f(x_0)\}$ is bounded.

Our method is based on monotonically reducing bounds on the error 
$f(x_b)-\wh f$ of the function value of the currently best point $x_b$. 
These bounds are derived from suitable linear relaxations and 
inequalities obtained with the help of a prox function. The solvability
of an auxiliary optimization subproblem involving the prox function is 
assumed. In many cases, this auxiliary subproblem has a cheap, closed 
form solution; this is shown here for the unconstrained case with a 
quadratic prox function, and in \sca{Ahookhosh \& Neumaier} 
\cite{AhoN1,AhoN2} for more general cases involving simple, practically
important convex sets $C$.

The OSGA algorithm presented here provides a fully adaptive alternative 
to current optimal first order methods. 
If no strong convexity is assumed, it shares the uniformity, the 
freeness of global parameters, and the optimal complexity properties 
of Lan's method, but has a far simpler structure and derivation. 
Beyond that, it also gives the optimal complexity in the 
strongly convex case, though it needs in this case -- like all other 
known methods with provable optimal linear convergence rate -- the 
knowledge of an explicit constant of strong convexity.
Furthermore -- like Nesterov's $O(k^{-2})$ algorithm from 
\cite{Nes.smooth} for the smooth case, but unlike his linearly 
convergent algorithm for the strongly convex case, scheme (2.2.19) in 
\sca{Nesterov} \cite{Nes.intL}, the algorithm derived here does not 
evaluate $f$ and $g$ outside their domain.
The method for analyzing the complexity of OSGA is also new; neither 
Tseng's complexity analysis nor Nesterov's estimating sequences are 
applicable to OSGA. 

The OSGA algorithm can be used in place of Nesterov's optimal algorithms
for smooth convex optimization and its variants whenever the latter are
traditionally employed. Thus it may be used as the smooth solver with 
methods for solving nonsmooth convex problems via smoothing 
(\sca{Nesterov} \cite{Nes.smooth}), and for solving very large linear 
programs (see, e.g., \sca{Aybat \& Iyengar} \cite{AybI},
\sca{Chen \& Burer} \cite{CheB}, \sca{Bu} et al.\ \cite{GuLW},
\sca{Nesterov} \cite{Nes.round,Nes.rel}, \sca{Richtarik} \cite{Ric.imp})

Numerical results are reported in \sca{Ahookhosh} \cite{Aho};
see also \sca{Ahookhosh \& Neumaier} \cite{AhoN.roks}.

\bigskip
{\bf Acknowledgment.} I'd like to thank Masoud Ahookhosh for numerous
useful remarks on an earlier version of the manuscript.

\section{The OSGA algorithm}

In this section we motivate and formulate the new algorithm.

In the following, $V$ denotes a Banach space with norm $\|\cdot\|$,
and $V^*$ is the dual Banach space with the dual norm $\|\cdot\|^*$.
$C$ is a closed, conves subset of $V$. The objective function 
$f:C\to \Rz$ is assumed to be convex, and $g(x)$ denotes a particular 
computable subgradient of $f$ at $x\in C$.

{\bf The basic idea.} 
The method is based on monotonically reducing bounds on the error 
$f(x_b)-\wh f$ of the function value of the currently best point $x_b$. 
These bounds are derived from suitable linear relaxations
\lbeq{e.f1}
f(z)\ge \gamma +\<h,z\> \Forall z\in C
\eeq
(where $\gamma\in\Rz$ and $h\in V^*$) with the 
help of a continuously differentiable \bfi{prox function} $Q:C\to \Rz$ 
satisfying 
\lbeq{e.Qinf}
Q_0:=\inf_{z\in C} Q(z) >0,
\eeq
\lbeq{e.strc}
Q(z)\ge Q(x)+\<g_Q(x),z-x\>+\frac{1}{2}\|z-x\|^2 \Forall x,z\in C,
\eeq
where $g_Q(x)$ denotes the gradient of $Q$ at $x\in C$. 
(Thus $Q$ is strongly convex with strong convexity parameter $\sigma=1$.
Choosing $\sigma=1$ simplified the formulas, and is no restriction of 
generality, as we may always rescale a prox function to enforce 
$\sigma=1$.)
We require that, for each $\gamma\in\Rz$ and $h\in V^*$, 
\lbeq{e.Eeta}
E(\gamma,h):=-\inf_{z\in C} \frac{\gamma+\<h,z\>}{Q(z)}
\eeq
is attained at some $z=U(\gamma,h)\in C$.
This requirement implies that, for arbitrary $\gamma_b\in\Rz$ and 
$h\in V^*$,
\[
\gamma_b+\<h,z\>\ge -E(\gamma_b,h) Q(z) \Forall z\in C.
\]
From \gzit{e.f1} for $z=\wh x$ and \gzit{e.Qinf}, we find that
\lbeq{e.ub}
\gamma_b=\gamma-f(x_b),~~~E(\gamma_b,h)\le \eta \implies 
0\le f(x_b) -\wh f\le \eta Q(\wh x).
\eeq
Typically,
\lbeq{e.u}
u:=U(\gamma_b,h),~~~ \ol\eta:=E(\gamma_b,h)
\eeq
are computed together; clearly one may update a given $\eta$ by 
$\ol\eta$, thus improving the bound. 

Note that the form of the auxiliary optimization problem \gzit{e.Eeta} 
is forced by this argument. Although this is a nonconvex optimization 
problem, it is shown in \sca{Ahookhosh \& Neumaier} \cite{AhoN1,AhoN2}
that there are many important cases where $E(\gamma,h)$ and 
$U(\gamma,h)$ are cheap to compute. In particular, we shall show in
Section \ref{s.prox} that this is the case when $C=V$ and the prox 
function is quadratic.

If an upper bound for $Q(\wh x)$ is known or assumed, the bound 
\gzit{e.ub} translates into a computable error estimate for the minimal 
function value. 
But even in the absence of such an upper bound, we can solve the 
optimization problem \gzit{e.f0} to a target accuracy
\lbeq{e.target}
0\le f(x_b) -\wh f\le\eps Q(\wh x)
\eeq
if we manage to decrease the error factor $\eta$ from its initial value 
until $\eta\le \eps$ for some target tolerance $\eps>0$.
This will be achieved by Algorithm \ref{a.OSGA} defined below. 
We shall prove for this algorithm complexity bounds on the number of 
iterations that are independent of the dimension of $V$ (which may be 
infinite), and -- apart from a constant factor -- best possible under 
a variety of assumptions on the objective function.

\bigskip
{\bf Constructing linear relaxations.}
The convexity of $f$ implies for $x,z\in C$ the bound
\lbeq{e.fzxb}
f(z)\ge f(x)+\<g(x),z-x\>,
\eeq
where $g(x)$ denotes a subgradient of $f$ at $x\in C$. Therefore 
\gzit{e.f1} always holds with
\[
\gamma=f(x_b)-\<g(x_b),x_b\>,~~~h=g(x_b).
\]
We can find more general relaxations of the form \gzit{e.f1} by 
accumulating past information. Indeed, if \gzit{e.f1} holds, 
$\alpha\in[0,1]$, and $x\in C$ then \gzit{e.f1} remains valid when we
substitute 
\[
\ol \gamma:= \gamma+\alpha(f(x)-\<g(x),x\>-\gamma),
\]
\[
\ol h := h+\alpha(g(x)-h)
\]
in place of $\gamma, h$, as by \gzit{e.fzxb},
\[
\bary{lll}
f(z)&=&(1-\alpha)f(z)+\alpha f(z)\\
&\ge& (1-\alpha)(\gamma+\<h,z\>)
+ \alpha(f(x)+\<g(x),z-x\>)\\
&=& (1-\alpha)\gamma+\alpha(f(x)-\<g(x),x\>) 
+\<(1-\alpha)h+\alpha g(x),z\>\\
&=& \ol \gamma+\<\ol h,z\>.
\eary
\]
For appropriate choices of $x$ and $\alpha$, this may give much 
improved error bounds. We discuss suitable choices for $x$ later.

\bigskip
{\bf Step size selection.}
The \bfi{step size parameter} $\alpha$ controls the fraction of the 
new information \gzit{e.fzxb} incorporated into the new relaxation.
It is chosen with the hope for a reduction factor of approximately 
$1-\alpha$ in the current error factor $\eta$, and must therefore be 
adapted to the actual progress made. 

First we note that in practice, $Q(\wh x)$ is unknown; hence the 
numerical value of $\eta$ is meaningless in itself. However, quotients 
of $\eta$ at different iterations have a meaning, quantifying the
amount of progress made.

In the following, we use bars to denote quantities tentatively 
modified in the current iteration, but they replace the current values 
of these quantities only if an acceptance criterion is met that we now 
motivate. We measure progress in terms of the quantity
\lbeq{e.R}
R:=\frac{\eta-\ol \eta}{\lambda\alpha \eta},
\eeq
where $\lambda\in{]0,1[}$ is a fixed number. A value $R\ge1$ indicates 
that we made sufficient progress in that 
\lbeq{e.oleR}
\ol \eta=(1-\lambda R\alpha)\eta
\eeq
was reduced at least by a fraction $\lambda$ of the designed improvement
of $\eta$ by $\alpha \eta$; thus the step size is acceptable or may 
even be increased if $R>1$. 
On the other hand, if $R<1$, the step size must be 
reduced significantly to improve the chance of reaching the design goal.
Introducing a maximal step size $\alpha_{\max}\in{]0,1[}$ and two 
parameters with $0<\kappa'\le\kappa$ to control the amount of increase 
or decrease in $\alpha$, we update the step size according to
\lbeq{e.olalp}
\ol\alpha:=\cases{\alpha e^{-\kappa} & if $R<1$,\cr
            \min(\alpha e^{\kappa' (R-1)},\alpha_{\max}) & if $R\ge1$.}
\eeq
Since updating the linear relaxation and $u$ makes sense only when 
$\eta$ was improved, we obtain the following update scheme.

\newpage
\begin{alg}\label{a.update} 
(\bfi{Update scheme})\\
global tuning parameters:~ $\lambda\in{]0,e^{-\kappa}]}$; 
                         ~ $\alpha_{\max}\in{]0,1[}$;
                         ~ $0<\kappa'\le\kappa$;\\
input:~ $\alpha,\eta,\ol h,\ol \gamma,\ol \eta,\ol u$;\\
output:~$\alpha,h,\gamma,\eta,u$;\\
$R=(\eta-\ol \eta)/(\lambda\alpha \eta)$;\\
if $R<1$, $\ol\alpha=\alpha e^{-\kappa}$; \\
else \spd $\ol\alpha=\min(\alpha e^{\kappa' (R-1)},\alpha_{\max})$;\\
end;\\
$\alpha= \ol \alpha$;\\
if $\ol \eta<\eta$,\\
\spc $h= \ol h$;~$\gamma= \ol \gamma$;
         ~$\eta= \ol \eta$;~$u= \ol u$;\\
end;\\
\end{alg}

\vspace*{-0.8cm}
If $\alpha_{\min}$ denotes the smallest actually occurring step size
(which is not known in advance), we have global linear convergence 
with a convergence factor of $1-e^{-\kappa}\alpha_{\min}$.
However, $\alpha_{\min}$ and hence this global rate of convergence 
may depend on the target tolerance $\eps$; thus the convergence speed 
in the limit $\eps\to 0$ may be linear or sublinear depending on the
properties of the specific function minimized.

\bigskip
{\bf Strongly convex relaxations.}
If $f$ is strongly convex, we may know a number $\mu>0$ such that
$f-\mu Q$ is still convex. In this case, we have in place of 
\gzit{e.fzxb} the stronger inequality
\lbeq{e.fzmu}
f(z)-\mu Q(z) \ge f(x)-\mu Q(x)+\<g(x)-\mu g_Q(x),z-x\> \for x,z\in C.
\eeq
In the following, we only assume that $\mu\ge 0$, thus covering the 
case of linear relaxations, too.

\gzit{e.fzmu} allows us to construct strongly convex relaxations of 
the form
\lbeq{e.fmu}
f(z)\ge \gamma +\<h,z\> +\mu Q(z) \Forall z\in C.
\eeq
For example, \gzit{e.fmu} always holds with
\[
h=g(x_b)-\mu g_Q(x_b),~~~\gamma=f(x_b)-\mu Q(x_b)-\<h,x_b\>.
\]
Again more general relaxations of the form \gzit{e.fmu} are found by 
accumulating past information. 

\begin{prop}\label{p.BGrad} 
Suppose that $x\in C$, $\alpha\in[0,1]$, and let
\[
\ol h = h+\alpha(g-h),~~~
\ol \gamma= \gamma
      +\alpha\Big(f(x)-\mu Q(x)-\<g,x\>-\gamma\Big),
\]
where
\[
g=g(x)-\mu g_Q(x).
\]
If \gzit{e.fmu} holds and $f-\mu Q$ is convex then \gzit{e.fmu} also 
holds with $\ol \gamma$ and $\ol h$ in place of $\gamma$ and $h$.
\end{prop}

\bepf
By \gzit{e.fzmu} and the assumptions,
\[
\bary{lll}
f(z)-\mu Q(z)&=&(1-\alpha)(f(z)-\mu Q(z))+\alpha(f(z)-\mu Q(z))\\
&\ge& (1-\alpha)(\gamma+\<h,z\>)\\
  &&   + \alpha\Big(f(x)-\mu Q(x)+\<g(x),z-x\>-\mu\<g_Q(x),z-x\>\Big)\\
&=& \ol \gamma+\<\ol h,z\>.
\eary
\]
\epf

The relaxations \gzit{e.fmu} lead to the following error bound.

\begin{prop}
Let
\[
\gamma_b:=\gamma-f(x_b),~~~
\eta:=E(\gamma_b,h)-\mu.
\]
Then \gzit{e.fmu} implies
\lbeq{e.bound}
0\le f(x_b) -\wh f\le \eta Q(\wh x).
\eeq

\end{prop}
\bepf
By definition of $E(\gamma_b,h)=\eta+\mu$ and \gzit{e.fmu}, we have 
\[
-(\eta+\mu)Q(z)\le \gamma_b+\<h,z\>=\gamma-f(x_b)+\<h,z\>
\le f(z)-f(x_b)-\mu Q(z).
\]
for all $z\in C$. Substituting $z=\wh x$ gives \gzit{e.bound}.
\epf

Note that for $\mu=0$, we simply recover the previous results for 
general convex functions.

\bigskip
{\bf An optimal subgradient algorithm.}
For a nonsmooth convex function, the subgradient at a point does not 
always determine a direction of descent. However, we may hope to find 
better points by moving from the best point $x_b$ into the direction 
of the point \gzit{e.u} used to determine our error bound. 
We formulate on this basis the following algorithm, for which optimal 
complexity bounds will be proved in Section \ref{s.boundsIter}.

\newpage
\begin{alg}\label{a.OSGA} 
(\bfi{Optimal subgradient algorithm}, \bfi{OSGA})\\
global tuning parameters:~ $\lambda,\alpha_{\max}\in{]0,1[}$;
                         ~ $0<\kappa'\le\kappa$;\\
input parameters: $\mu\ge 0$;~$\eps>0$;~$f_\target$;\\
output:~$x_b$;\\
assumptions: $f-\mu Q$ is convex;\\
begin\\
\spc choose $x_b$;~stop if $f(x_b)\le f_\target$;\\
\spc $h=g(x_b)-\mu g_Q(x_b)$;~$\gamma=f(x_b)-\mu Q(x_b)-\<h,x_b\>$;\\
\spc $\gamma_b=\gamma-f(x_b)$;
     $u=U(\gamma_b,h)$;~$\eta=E(\gamma_b,h)-\mu$;\\
\spc $\alpha=\alpha_{\max}$;\\
\spc while 1,\\
\spc \spc $x=x_b+\alpha(u-x_b)$;~$g=g(x)-\mu g_Q(x)$;\\
\spc \spc $\ol h = h+\alpha(g-h)$;
          $\ol\gamma=\gamma+\alpha(f(x)-\mu Q(x)-\<g,x\>-\gamma)$;\\
\spc \spc $x_b' =\D\argmin_{z\in\{x_b,x\}} f(z)$;\\
\spc \spc $\gamma_b'=\ol\gamma-f(x_b')$;~
          $u'=U(\gamma_b',\ol h)$;~$x'=x_b+\alpha(u'-x_b)$;\\
\spc \spc choose $\ol x_b$ with $f(\ol x_b)\le \min(f(x_b'),f(x'))$;\\
\spc \spc $\ol\gamma_b=\ol\gamma-f(\ol x_b)$;~
          $\ol u=U(\ol \gamma_b,\ol h)$;~
          $\ol \eta= E(\ol\gamma_b,\ol h)-\mu$;\\
\spc \spc $x_b= \ol x_b$;\\
\spc \spc stop if some user-defined test is passed;\\
\spc \spc update $\alpha$, $h$, $\gamma$, $\eta$, $u$ 
          by Algorithm \ref{a.update};\\
\spc end;\\
end;
\end{alg}

Note that the strong convexity parameter $\mu$ needs to be specified to 
use the algorithm. If $\mu$ is unknown, one may always put $\mu=0$ 
(ignoring possible strong convexity), at the cost of possibly slower 
worst case asymptotic convergence. 
(Techniques like those used in \sca{Juditsky \& Nesterov} \cite{JudN} 
or \sca{Gonzaga \& Karas} \cite{GonK} for choosing $\mu$ adaptively can 
probably be applied to the above algorithm to remove the dependence on 
having to know $\mu$. However, \cite{JudN} requires an explicit 
knowledge of a Lipschitz constant for the gradient, while \cite{GonK} 
proves only sublinear convergence. It is not yet clear how to avoid 
both problems.)

The analysis of the algorithm will be independent of the choice of 
$\ol x_b$ allowed in Algorithm \ref{a.OSGA}. The simplest choice is 
$\ol x_b =\D\argmin_{z\in\{x_b',x'\}} f(z)$. If the best function 
value $f(x_b)$ is stored and updated, each iteration then requires the 
computation of two function values $f(x)$ and $f(x')$ and one 
subgradient $g(x)$. 

However, the algorithm allows the incorporation of heuristics 
to look for improved function values before deciding on the choice of 
$\ol x_b$. This may involve additional function 
evaluations at points selected by a line search procedure (see, e.g., 
\sca{Beck \& Teboulle} \cite{BecT}), a bundle optimization (see, e.g., 
\sca{Lan} \cite{Lan}), or a local quadratic approximation (see, e.g.,
\sca{Yu} et al. \cite{YuVGS}).

Numerical results are reported in \sca{Ahookhosh} \cite{Aho};
see also \sca{Ahookhosh \& Neumaier} \cite{AhoN.roks}.

\section{Inequalities for the error factor}

The possibility to get worst case complexity bounds rests on the
establishment of a strong upper bound on the error factor $\eta$.
This bound depends on global information about the function $f$;
while not necessary for executing the algorithm itself, it is needed
for the analysis. Depending on the properties of $f$,
global information of different strength can be used, resulting in 
inequalities of corresponding strength. The key steps in the analysis
rely on the following lower bound for the term $\gamma+\<h,z\>$.

\begin{prop}\label{p.EU}
Let $v=U(\gamma,h)$. Then 
\lbeq{e.etae}
\gamma+\<h,v\>=-E(\gamma,h)Q(v).
\eeq
Moreover, if $E(\gamma,h)\ge 0$ then for all $z\in C$, 
\lbeq{e.key}
\gamma+\<h,z\> \ge E(\gamma,h)\Big(\half\|z-v\|^2-Q(z)\Big),
\eeq
\lbeq{e.key2}
 E(\gamma,h)(Q(z)-Q(v))+\<h,z-v\> \ge 0.
\eeq

\end{prop}

\bepf
By definition of $E(\gamma,h)$, the function $\phi:C\to \Rz$ defined by
\[
\phi(z):=\gamma+\<h,z\>+E(\gamma,h)Q(z)
\]
is nonnegative and vanishes for $z=v:=U(\gamma,h)$. This implies 
\gzit{e.etae}. Writing $g_\phi(z)$ for the gradient of $\phi$ at
$z\in C$, strong convexity \gzit{e.strc} of $Q$ implies
\[
\phi(z)-\phi(v)-\<g_\phi(v),z-v\>
=E(\gamma,h)\Big(Q(z)-Q(v)-\<g_Q(v),z-v\>\Big)
\ge \frac{E(\gamma,h)}{2}\|z-v\|^2.
\]
But 
\[
\<g_\phi(v),z-v\>
=\lim_{\alpha\downto 0} \frac{\phi(v+\alpha(z-v))-\phi(v)}{\alpha} \ge 0
\]
since $\phi(v)=0$. This proves \gzit{e.key}. If we eliminate $\gamma$
using \gzit{e.etae} and delete the norm term, we obtain \gzit{e.key2}. 
\epf

\begin{thm}\label{t.tol1}
In Algorithm \ref{a.OSGA}, the error factors are related by 
\lbeq{e.oleN}
\ol \eta -(1-\alpha)\eta
\le \frac{\alpha^2\|g(x)\|_*^2}{2(1-\alpha)(\eta+\mu) Q_0},
\eeq
where $\|\cdot\|_*$ denotes the norm dual to $\|\cdot\|$.
\end{thm}

\bepf
We first establish some inequalities needed for the later 
estimation.
By convexity of $Q$ and the definition of $\ol h$,
\[
\bary{lll}
\alpha\mu\Big(Q(\ol u)- Q(x)+\<g_Q(x),x\>\Big)
&\ge&  \alpha\mu\<g_Q(x),\ol u\>
= \<h-\ol h+\alpha(g(x)-h),\ol u\>\\
&=&(1-\alpha)\<h,\ol u\>+\<\alpha g(x)-\ol h,\ol u\>.
\eary
\]
By definition of $x$, we have
\[
(1-\alpha)(x_b-x)=-\alpha(u-x).
\]
Hence \gzit{e.fzmu} (with $\mu=0$) implies
\[
(1-\alpha)(f(x_b)-f(x))\ge(1-\alpha)\<g(x),x_b-x\>=-\alpha\<g(x),u-x\>.
\]
By definition of $\ol \gamma$, we conclude from these two inequalities 
that
\[
\bary{lll}
\ol \gamma-f(x)+\alpha\mu Q(\ol u)
&=& (1-\alpha)(\gamma-f(x))-\alpha\<g(x),x\>
     +\alpha\mu\Big(Q(\ol u)- Q(x)+\<g_Q(x),x\>\Big)\\
&\ge& (1-\alpha)\Big(\gamma-f(x)+\<h,\ol u\>\Big)
       +\alpha\<g(x),\ol u -x\> - \<\ol h,\ol u\>\\
&\ge& (1-\alpha)\Big(\gamma-f(x_b)+\<h,\ol u\>\Big)
      +\alpha\<g(x),\ol u -u\> -\<\ol h,\ol u\>.\\
\eary
\]
Using this, \gzit{e.etae} (with $\ol\gamma_b=\ol\gamma-f(\ol x_b)$ in 
place of $\gamma$ and $\ol h$ in place of $h$), and 
$E(\ol \gamma_b, \ol h)=\ol \eta+\mu$ now gives
\lbeq{e.f14}
\bary{lll}
(\ol \eta+\mu-\alpha\mu)Q(\ol u)
&=&f(\ol x_b)-\ol\gamma-\<\ol h,\ol u\>-\alpha\mu Q(\ol u)\\
&\le&f(\ol x_b)-f(x)-\alpha\<g(x),\ol u-u\>
    -(1-\alpha)\Big(\gamma-f(x_b)+\<h,\ol u\>\Big).
\eary
\eeq
Using \gzit{e.key} (with $\gamma_b=\gamma-f(x_b)$ in place of $\gamma$) 
and $\eta+\mu=E(\gamma_b,h)$, we find
\lbeq{e.emq}
(\eta+\mu)Q(\ol u)
\ge f(x_b)-\gamma-\<h,\ol u\> +\frac{\eta+\mu}{2}\|\ol u-u\|^2.
\eeq
Now \gzit{e.f14} and \gzit{e.emq} imply
\[
\bary{lll}
(\ol\eta -(1-\alpha)\eta)Q(\ol u)
&=& (\ol\eta +\mu-\alpha\mu)Q(\ol u) -(1-\alpha)(\eta+\mu)Q(\ol u)\\
&\le& f(\ol x_b)-f(x)-(1-\alpha)\Big(\gamma-f(x_b)+\<h,\ol u\>\Big)\\
&&    -\alpha\<g(x),\ol u-u\> \\
&&    -(1-\alpha)\Big(
   f(x_b)-\gamma-\<h,\ol u\>+\D\frac{\eta+\mu}{2}\|\ol u-u\|^2\Big)
\\[3mm]
&=& f(\ol x_b)-f(x)+\ol S,
\eary
\]
where
\lbeq{e.Delta}
\bary{ll}
\ol S&:=-\alpha\<g(x),\ol u-u\>
   -\D\frac{(1-\alpha)(\eta+\mu)}{2}\|\ol u-u\|^2\\[4mm]
 &\le\alpha\|g(x)\|_*\|\ol u-u\|
   -\D\frac{(1-\alpha)(\eta+\mu)}{2}\|\ol u-u\|^2\\[4mm]
&= \D\frac{\alpha^2 \|g(x)\|_*^2
           -(\alpha \|g(x)\|_*+(1-\alpha)(\eta+\mu)\|\ol u-u\|)^2}
           {2(1-\alpha)(\eta+\mu)}
\le \D\frac{\alpha^2\|g(x)\|_*^2}{2(1-\alpha)(\eta+\mu)}.
\eary
\eeq
If $\ol\eta \le(1-\alpha)\eta$ then \gzit{e.oleN} holds trivially.
Thus we assume that $\ol\eta >(1-\alpha)\eta$. Then 
\lbeq{e.taud}
(\ol\eta -(1-\alpha)\eta)Q_0  
\le (\ol\eta -(1-\alpha)\eta)Q(\ol u) 
\le  f(\ol x_b)-f(x)+\ol S.
\eeq
Since $f(\ol x_b)\le f(x)$, we conclude again that
\gzit{e.oleN} holds. Thus \gzit{e.oleN} holds generally.
\epf

Note that the arguments used in this proof did not make use of $x'$; 
thus \gzit{e.oleN} even holds when one sets $x'=x$ in the algorithm, 
saving some work.

\begin{thm}\label{t.tol2}
If $f$ has Lipschitz continuous gradients with Lipschitz constant $L$ 
then, in Algorithm \ref{a.OSGA},
\lbeq{e.oleS}
 \ol\eta > (1-\alpha)\eta 
\implies (1-\alpha)(\eta+\mu) < \alpha^2L.
\eeq
\end{thm}

\bepf
The proof follows the general line of the preceding proof, but now
we must consider the information provided by $x'$.

Since $E$ is monotone decreasing in its first argument
and $f(x_b')\ge f(\ol x_b)$, the hypothesis of \gzit{e.oleS} implies 
that
\[
\eta':= E(\ol \gamma-f(x_b'),\ol h)-\mu 
\ge E(\ol \gamma-f(\ol x_b),\ol h)-\mu 
     = \ol\eta > (1-\alpha)\eta.
\]
By convexity of $Q$ and the definition of $\ol h$,
\[
\bary{lll}
\alpha\mu\Big(Q(u')- Q(x)+\<g_Q(x),x\>\Big)
&\ge&  \alpha\mu\<g_Q(x),u'\>
= \<h-\ol h+\alpha(g(x)-h),u'\>\\
&=&(1-\alpha)\<h,u'\>+\<\alpha g(x)-\ol h,u'\>.
\eary
\]
By definition of $x$, we have
\[
(1-\alpha)(x_b-x)=-\alpha(u-x).
\]
Hence \gzit{e.fzmu} (with $\mu=0$) implies
\[
(1-\alpha)(f(x_b)-f(x))\ge(1-\alpha)\<g(x),x_b-x\>=-\alpha\<g(x),u-x\>.
\]
By definition of $\ol \gamma$, we conclude from the last two 
inequalities that
\[
\bary{lll}
\ol \gamma-f(x)+\alpha\mu Q(u')
&=& (1-\alpha)(\gamma-f(x))-\alpha\<g(x),x\>
     +\alpha\mu\Big(Q(u')- Q(x)+\<g_Q(x),x\>\Big)\\
&\ge& (1-\alpha)\Big(\gamma-f(x)+\<h,u'\>\Big)
       +\alpha\<g(x),u' -x\> - \<\ol h,u'\>\\
&\ge& (1-\alpha)\Big(\gamma-f(x_b)+\<h,u'\>\Big)
      +\alpha\<g(x),u' -u\> -\<\ol h,u'\>.\\
\eary
\]
Using this, \gzit{e.etae} (with $\gamma_b'=\ol\gamma-f(x_b')$ in 
place of $\gamma$ and $\ol h$ in place of $h$), and 
$E(\gamma_b', \ol h)=\eta'+\mu$ now gives
\lbeq{e.f142}
\bary{lll}
(\eta'+\mu-\alpha\mu)Q(u')
&=&f(x_b')-\ol\gamma-\<\ol h,u'\>-\alpha\mu Q(u')\\
&\le&f(x_b')-f(x)-\alpha\<g(x),u'-u\>\\
&&    -(1-\alpha)\Big(\gamma-f(x_b')+\<h,u'\>\Big).
\eary
\eeq
Using \gzit{e.key} (with $\gamma_b=\gamma-f(x_b)$ in place of $\gamma$) 
and $\eta+\mu=E(\gamma_b,h)$, we find
\lbeq{e.emq2}
(\eta+\mu)Q(u')
\ge f(x_b)-\gamma-\<h,u'\> +\frac{\eta+\mu}{2}\|u'-u\|^2.
\eeq
Now \gzit{e.f142} and \gzit{e.emq2} imply
\[
\bary{lll}
(\eta' -(1-\alpha)\eta)Q(u')
&=& (\eta' +\mu-\alpha\mu)Q(u') -(1-\alpha)(\eta+\mu)Q(u')\\
&\le& f(x_b')-f(x)-(1-\alpha)\Big(\gamma-f(x_b)+\<h,u'\>\Big)\\
&&    -\alpha\<g(x),u'-u\> \\
&&    -(1-\alpha)\Big(
   f(x_b)-\gamma-\<h,u'\>+\D\frac{\eta+\mu}{2}\|u'-u\|^2\Big)
\\[3mm]
&=& f(x_b')-f(x)+S',
\eary
\]
where
\[
S':=-\alpha\<g(x),u'-u\>-\frac{(1-\alpha)(\eta+\mu)}{2}\|u'-u\|^2,
\]
giving
\[
(\eta' -(1-\alpha)\eta)Q_0  \le  f(x_b')-f(x)+S'.
\]
Now 
\lbeq{e.Lipgrad}
\bary{lll}
f(x_b')&\le& \D f(x')\le f(x)+\<g(x),x'-x\>+\frac{L}{2}\|x'-x\|^2\\
&=&\D f(x)+\alpha\<g(x),u'-u\>+\frac{\alpha^2 L}{2}\|u'-u\|^2,
\eary
\eeq
so that under the hypothesis of \gzit{e.oleS}
\[
0 < (\eta' -(1-\alpha)\eta)Q_0
\le \frac{\alpha^2L-(1-\alpha)(\eta+\mu)}{2}\|u'-u\|^2.
\]
Thus $\alpha^2L-(1-\alpha)(\eta+\mu)>0$, and the conclusion of 
\gzit{e.oleS} holds.
\epf

\section{Bounds for the number of iterations}\label{s.boundsIter}

We now use the inequalities from Theorem \ref{t.tol1} and Theorem 
\ref{t.tol2} to derive bounds for the number of iterations. 
The weakest global assumption, mere convexity, leads to the weakest 
bounds and guarantees sublinear convergence only, while the strongest 
global assumption, strong convexity and Lipschitz continuous gradients, 
leads to the strongest bounds guaranteeing $R$-linear convergence. 
Our main result shows that, asymptotically as $\eps\to 0$,
the number of iterations needed by the OSGA algorithm matches the 
lower bounds on the complexity derived by \sca{Nemirovski \& Yudin} 
\cite{NemY}, apart from constant factors:

\begin{thm}\label{t.cBSGA}
Suppose that $f-\mu Q$ is convex. Then:

(i) (\bfi{Nonsmooth complexity bound})\\
If the points generated by Algorithm \ref{a.OSGA} stay in a bounded 
region of the interior of $C$, or if $f$ is Lipschitz continuous in $C$,
the total number of iterations needed to reach a point with 
$f(x)\le f(\wh x)+\eps$ is at most $O((\eps^2+\mu\eps)^{-1})$. 
Thus the asymptotic worst case complexity is $O(\eps^{-2})$ when 
$\mu=0$ and $O(\eps^{-1})$ when $\mu>0$.

(ii) (\bfi{Smooth complexity bound})\\
If $f$ has Lipschitz continuous gradients with Lipschitz constant $L$, 
the total number of iterations needed by Algorithm \ref{a.OSGA} to 
reach a point with $f(x)\le f(\wh x)+\eps$ is at most $O(\eps^{-1/2})$ 
if $\mu=0$, and at most $\D O(|\log\eps|\sqrt{L/\mu})$ if $\mu>0$.
\end{thm}

In particular, if $f$ is strongly convex and differentiable with 
Lipschitz continuous gradients, $\mu>0$ 
holds with arbitrary quadratic prox functions, and we get a complexity 
bound similar to that achieved by the preconditioned conjugate gradient 
method for linear systems; cf. \sca{Axelsson \& Lindskog} \cite{AxeL}.

Note that \gzit{e.oleS} generalizes to other situations by replacing 
\gzit{e.Lipgrad} with a weaker smoothness property of the form 
\lbeq{e.weakgrad}
f(z) \le f(x)+\<g(x),z-x\>+ \phi(\|z-x\|)
\eeq
with $\phi$ convex and monotone increasing. For example, this holds 
with $\phi(t)=L_1t$ if $f$ has subgradients with bounded variation, 
and with $\phi(t)=L_s t^{s+1}$ if $f$ has H\"older continuous gradients 
with exponent $s\in{]0,1[}\,$, and with linear combinations thereof in 
the composite case considered by \sca{Lan} \cite{Lan}. Imitating
the analysis below of the two cases stated in the theorem then gives 
corresponding complexity bounds matching those obtained by Lan.

Theorem \ref{t.cBSGA} follows from the two propositions below
covering the different cases, giving in each case explicit upper bounds 
on the number $K_\mu(\alpha,\eta)$ of further iterations needed to 
complete the algorithm form a point where the values of $\alpha$ and 
$\eta$ given as arguments of $K_\mu$ were achieved. 
We write $\alpha_0$ and $\eta_0$ for the initial values of $\alpha$ and 
$\eta$. 
Only the dependence on $\mu$, $\alpha$, and $\eta$ is made explicit.

\begin{prop} \label{p.boundN}
Suppose that the dual norm of the subgradients $g(x)$ encountered 
during the iteration remains bounded by the constant $c_0$.
Let $c_1>0$, and define
\[
c_1:=\frac{c_0^2}{2Q_0},~~~
c_2:=\max\Big(\frac{c_1}{(1-e^{-\kappa})(1-\alpha_{\max})},
              \frac{\eta_0(\eta_0+\mu)}{\alpha_0}\Big),~~~
c_3=\frac{c_2}{2\lambda}.
\]
(i) In each iteration,
\lbeq{e.ebound0}
\eta(\eta+\mu)\le\alpha c_2.
\eeq
(ii) The algorithm stops after at most
\lbeq{e.kbound0}
K_\mu(\alpha,\eta)
:= 1 + \kappa^{-1}\log \frac{c_2\alpha}{\eps(\eps+\mu)}
     + \frac{c_3}{\eps(\eps+\mu)}
     -\frac{c_3}{\eta(\eta+\mu)}
\eeq
further iterations.

In particular, (i) and (ii) hold when the iterates stay in a bounded 
region of the interior of $C$, or when $f$ is Lipschitz continuous in 
$C$. 
\end{prop}

Note that any convex function is Lipschitz continuous in any closed and 
bounded domain inside its support. Hence if the iterates stay in a 
bounded region $R$ of the interior of $C$, $\|g\|$ is bounded by the 
Lipschitz constant of $f$ in the closure of the region $R$. 

\bepf
(i) 
Condition \gzit{e.ebound0} holds initially, and is preserved in each 
update unless $\alpha$ is reduced. But then $R<1$, hence
$\ol\eta \ge (1-\lambda\alpha)\eta$. Thus Theorem \ref{t.tol1} implies 
\[
(1-\lambda)\alpha \eta 
\le \ol\eta -(1-\alpha)\eta 
\le \frac{\alpha^2 c_1}{(1-\alpha)(\eta+\mu)}.
\]
This implies
\[
(1-\lambda)(1-\alpha) \eta(\eta+\mu) \le \alpha c_1,
\]
and since $\lambda\le e^{-\kappa}<1$,
\[
\ol\eta(\ol \eta+\mu)\le \eta(\eta+\mu) 
\le \frac{\alpha c_1}{(1-\lambda)(1-\alpha)}
\le\alpha c_2.
\] 
Thus \gzit{e.ebound0} holds again after the reduction, and hence always.

(ii) 
As the algorithm stops once $\eta\le \eps$, \gzit{e.ebound0} implies 
that in each iteration $c_2\alpha\ge\eps(\eps+\mu)$. As $\alpha$ is 
reduced only when $R<1$, and then by a fixed factor $e^{-\kappa}$, 
this cannot happen more than 
$\D\kappa^{-1}\log\frac{c_2\alpha}{\eps(\eps+\mu)}$ 
times in turn. Thus after some number of $\alpha$-reductions we must 
always have another step with $R\ge1$. 
By \gzit{e.oleR}, this gives a reduction of $\eta$ by a factor of at 
least $1-\lambda\alpha$. But this implies 
that the stopping criterion $\eta\le \eps$ is eventually reached.
Therefore the algorithm stops eventually.
Since $R\ge 0$, \gzit{e.olalp} implies 
$\ol \alpha \le \alpha e^{\kappa(R-1)}$. Therefore
\lbeq{e.alpbar}
\log(\alpha/\ol \alpha) \ge \kappa(1-R).
\eeq
Now \gzit{e.alpbar}, \gzit{e.oleR}, and \gzit{e.ebound0} imply
\[
\bary{lll}
K_\mu(\alpha,\eta)-K_\mu(\ol \alpha,\ol\eta)
&=& \D\frac{\log (\alpha/\ol\alpha)}{\kappa}
   +\frac{c_3}{\ol\eta(\ol \eta+\mu)}-\frac{c_3}{\eta(\eta+\mu)}\\[4mm]
&\ge& 
1-R+\D\frac{c_3}{\ol\eta(\ol \eta+\mu)}-\frac{c_3}{\eta(\eta+\mu)}
\\[4mm]
&=& 1-R
   +\D\frac{c_3}{(1-\lambda R\alpha)\eta((1-\lambda R\alpha)\eta+\mu)}
   -\frac{c_3}{\eta(\eta+\mu)}\\[4mm]
&\ge& 1-R +\D 
c_3\frac{(\eta+\mu)-(1-\lambda R\alpha)((1-\lambda R\alpha)\eta+\mu)}
{(1-\lambda R\alpha)\eta(\eta+\mu)((1-\lambda R\alpha)\eta+\mu)}
\\[4mm]
&=& 1-R
+\D c_3\frac{\lambda R\alpha((2-\lambda R\alpha)\eta+\mu)}
{(1-\lambda R\alpha)\eta(\eta+\mu)((1-\lambda R\alpha)\eta+\mu)}
\\[4mm]
&\ge& 1-R+\D c_3\frac{2\lambda R\alpha} {\eta(\eta+\mu)}
= 1-R+\D\frac{2c_3\lambda R}{c_2}=1.
\eary
\]
This implies the complexity bound by reverse induction, since 
immediately before the last iteration, $c_2\alpha\ge\eps(\eps+\mu)$ and
$\eta>\eps$, hence $K_\mu(\alpha,\eta)\ge 1$.
\epf

\begin{prop}\label{p.boundS}
Suppose that $f$ has Lipschitz continuous gradients with Lipschitz 
constant $L$, and put
\[
c_4=\max\Big(\frac{\eta_0+\mu}{\alpha_0^2},
             \frac{e^{2\kappa} L}{1-\alpha_{\max}}\Big),~~~
c_5=\frac{4c_4}{\lambda^2},~~~
c_6=\sqrt{\frac{c_4}{\mu}},~~~c_7=\frac{c_6}{\lambda}.
\]
(i) In each iteration
\lbeq{e.ebound1}
\eta+\mu\le \alpha^2c_4,
\eeq
(ii) The algorithm stops after at most $K_\mu(\alpha,\eta)$ further 
iterations. Here
\lbeq{e.kbound1}
K_0(\alpha,\eta)
:=1+\kappa^{-1}\log \Big(\alpha\sqrt{\frac{c_4}{\eps}}\Big)
   +\sqrt{\frac{c_5}{\eps}}-\sqrt{\frac{c_5}{\eta}},
\eeq
\lbeq{e.kbound2}
K_\mu(\alpha,\eta)
:=1+\frac{\log (c_6\alpha)}{\kappa}+c_7\log\frac{\eta}{\eps} \for \mu>0.
\eeq
\end{prop}

\bepf
(i) \gzit{e.ebound1} holds initially, and is preserved in each update 
unless $\alpha$ is reduced. But then $R<1$, hence by 
Theorem \ref{t.tol2}, $(1-\alpha)(\eta+\mu)= \alpha^2 L$ before 
the reduction. Therefore
\[
\ol\eta +\mu
\le \eta +\mu
\le \frac{\alpha^2 L}{1-\alpha} 
\le \frac{\alpha^2 L}{1-\alpha_{\max}}
\le \alpha^2e^{-2\kappa} c_4
\le \ol \alpha^2 c_4.
\] 
Thus \gzit{e.ebound1} holds again after the reduction, and hence always.
As in the previous proof, we find that the algorithm stops 
eventually, and \gzit{e.alpbar} holds.

(ii) If $\mu=0$ then \gzit{e.oleR}, \gzit{e.alpbar}, and 
\gzit{e.ebound1} imply
\[
\bary{lll}
K_0(\alpha,\eta)-K_0(\ol \alpha,\ol\eta)
&=& \D\frac{\log (\alpha/\ol\alpha)}{\kappa}
  +\sqrt{\frac{c_5}{(1-\lambda R\alpha)\eta}}
  -\sqrt{\frac{c_5}{\eta}}\\[4mm]
&\ge& 1-R+\D\Big(1-\sqrt{1-\lambda R\alpha}\Big)
        \sqrt{\frac{c_5}{(1-\lambda R\alpha)\eta}}\\[4mm]
&\ge& 1-R+\D\frac{\lambda R\alpha}{2}\sqrt{\frac{c_5}{\eta}}
\ge 1-R+\D\frac{\lambda R}{2}\sqrt{\frac{c_5}{c_4}}=1.
\eary
\]
This implies the complexity bound by reverse induction, since 
immediately before the last iteration, $\alpha\ge \sqrt{\eps/c_4}$, 
hence $K_0(\alpha,\eta)\ge 1$.

(iii) If $\mu>0$ then \gzit{e.ebound1} shows that always 
$c_6\alpha\ge 1$, hence 
$\ol\eta=(1-\lambda R\alpha)\eta \le (1-R/c_7)\eta$.
Therefore 
\[
K_\mu(\alpha,\eta)-K_\mu(\ol \alpha,\ol\eta)
= \frac{\log (\alpha/\ol\alpha)}{\kappa}+ c_7\log\frac{\eta}{\ol\eta}
\ge 1-R+ c_7\log\frac{1}{1-R/c_7} \ge 1,
\]
and the result follows as before.
\epf

{\em Proof} of Theorem \ref{t.cBSGA}.\\
(i) We apply Proposition \ref{p.boundN}(ii) to the first iteration,
and note that $K_0(\alpha,\eta)=O(e^{-2})$ and
$K_\mu(\alpha,\eta)=O(e^{-1})$ if $\mu>0$.\\
(i) We apply Proposition \ref{p.boundS}(ii) to the first iteration,
and note that $K_0(\alpha,\eta)=O(e^{-{1/2}})$ and
$K_\mu(\alpha,\eta)=O(\log \eps^{-1})$ if $\mu>0$.

\section{Quadratic prox functions for unconstrained problems}
\label{s.prox}

To use the algorithm in practice, we need prox functions for which
$E(\gamma,h)$ and $U(\gamma,h)$ can be evaluated easily. We first 
derive the optimality conditions for the asssociated auxiliary 
optimization problem \gzit{e.Eeta}.

\begin{prop}\label{p.E}
If the function $E_{\gamma,h}:C\to \Rz$ defined by
\[
E_{\gamma,h}(z):=-\frac{\gamma+\<h,z\>}{Q(z)}
\]
attains its supremum $E(\eta,h)$ at $z=u$, the point $u=U(\gamma,h)$ 
satisfies
\lbeq{e.etae2}
E(\gamma,h)Q(u)=-\gamma-\<h,u\>,
\eeq
\lbeq{e.ye}
\<E(\gamma,h)g_Q(u)+h,z-u\>\ge 0 \Forall z\in C.
\eeq
\end{prop}

\bepf
\gzit{e.etae2} holds since $E(\gamma,h)=E_{\gamma,h}(u)$.
To get \gzit{e.ye}, we differentiate the identity
\[
E_{\gamma,h}(z)Q(z)=-\gamma-\<h,z\>
\]
and obtain
\[
\frac{\partial E_{\gamma,h}(z)}{\partial z}Q(z)+E_{\gamma,h}(z)g_Q(z)
=-h.
\]
Now \gzit{e.ye} follows from $Q(z)>0$ and the first order optimality 
condition
\[
\Big\<\frac{-\partial E_{\gamma,h}(u)}{\partial z},z-u\Big\>\ge 0
 \Forall z\in C
\]
for a minimum of $-E_{\gamma,h}(z)$ at $z=u$.
\epf

Simple domains $C$ and prox functions for which the conditions
\gzit{e.etae2} and \gzit{e.ye} can be solved easily are discussed in
 \sca{Ahookhosh \& Neumaier} \cite{AhoN1,AhoN2}.

In the remainder, we only discuss the simplest case, where the original 
optimization problem is unconstrained (so that $C=V$) and the norm on 
$V$ is Euclidean,  
\[
\|z\|:=\sqrt{\<Bz,z\>},
\]
where the \bfi{preconditioner} $B$ is a symmetric and positive 
definite linear mapping $B:V\to V^*$.
The associated dual norm on $V^*$ is then given by 
\[
\|h\|_*:=\|B^{-1}h\|=\sqrt{\<h,B^{-1}h\>}.
\]
Given the preconditioner, it is natural to consider the quadratic
prox function
\lbeq{e.Q0}
Q(z):= Q_0+\half\|z-z_0\|^2,
\eeq
where $Q_0$ is a positive number and $z_0\in V$. 
The assumption of Proposition \ref{p.E} is satisfied since as the 
quotient of a linear and a positive quadratic function, 
$E_{\gamma,h}(z)$ takes positive and negative values and is arbitrarily 
small outside a ball of sufficiently large radius. Therefore the 
level sets for nonzero function values are compact, and the supremum 
of a continuous function on a compact set is attained. Similarly, the 
infimum is attained.

Since $C=V$ and $g_Q(z)=B(z-z_0)$, we conclude from the proposition that
$E(\gamma,h) B(u-z_0)+h=0$, where $u=U(\gamma,h)$, so that
\lbeq{e.h0}
U(\gamma,h)=z_0-E(\gamma,h)^{-1}B^{-1}h.
\eeq
Inserting this into \gzit{e.etae2} and writing $e=E(\gamma,h)$, we find 
\[
e\Big(Q_0+\half\|-e^{-1}B^{-1}h\|^2\Big)
=eQ(u)=-\gamma-\<h,z_0-e^{-1}B^{-1}h\>,
\]
which simplifies to the quadratic equation
\[
Q_0 e^2+\beta e-\half\|h\|_*^2=0,~~~\beta=\gamma+\<h,z_0\>.
\]
The two solutions are the only stationary points, hence the solution 
with positive (negative) function value must be the unique maximizer
(resp. minimizer). One easily checks that the maximizer is given by
\lbeq{e.E0}
E(\gamma,h)
=\frac{-\beta+\sqrt{\beta^2+2Q_0\|h\|_*^2}}{2Q_0}
=\frac{\|h\|_*^2}{\beta+\sqrt{\beta^2+2Q_0\|h\|_*^2}}.
\eeq
(The first form is numerically stable when $\beta\le0$, the second
when $\beta>0$.)

A reasonable choice is to take the starting point of the iteration for 
$z_0$, and $Q_0\approx\half\|\wh x-z_0\|^2$, an order of magnitude 
guess.




\begin{thebibliography}{99}

\bibitem{Aho} Ahookhosh, M.:
Optimal subgradient algorithms with application to large-scale linear
inverse problems, Manuscript, University of Vienna (2014).

\bibitem{AhoN.roks} M. Ahookhosh and A. Neumaier,
High-dimensional convex optimization via optimal affine subgradient 
algorithms, in ROKS workshop (2013), 83--84.

\bibitem{AhoN1} M. Ahookhosh and A. Neumaier,
Optimal subgradient methods for structured convex constrained 
optimization. I: theoretical results,
Manuscript, University of Vienna (2014).

\bibitem{AhoN2} M. Ahookhosh and A. Neumaier,
Optimal subgradient methods for structured convex constrained 
optimization. II: numerical results,
Manuscript, University of Vienna (2014).

\bibitem{AusT} A. Auslender and M. Teboulle,
Interior gradient and proximal methods for convex and conic 
optimization,
SIAM J. Optimization 16 (2006), 697--725.

\bibitem{AxeL} O. Axelsson and G. Lindskog, 
On the rate of convergence of the conjugate gradient method,
Numer. Math. 48 (1986), 499--523.

\bibitem{AybI} N.S. Aybat and G. Iyengar,
A first-order augmented Lagrangian method for compressed sensing,
SIAM J. Optim., 22(2) (2012), 429â459.

\bibitem{BecT} A. Beck and M. Teboulle,
A fast iterative shrinkage-thresholding algorithm for linear inverse 
problems,
SIAM J. Imaging Sci. 2 (2009), 183--202.

\bibitem{BecCG} S.R. Becker, E.J. Cand{\`e}s, and M.C. Grant, 
Templates for convex cone problems with applications to sparse signal 
recovery,
Math. Programming Comput. 3 (2011), 165--218.

\bibitem{CheB} J. Chen and S. Burer,
A First-Order Smoothing Technique for a Class of Large-Scale Linear 
Programs,
Manuscript (2011).
\url{http://www.optimization-online.org/DB_FILE/2011/11/3233.pdf}

\bibitem{DevGN} O. Devolder, F. Glineur, and Y. Nesterov,
First-order methods of smooth convex optimization with inexact oracle,
Mathematical Programming, (2013) DOI 10.1007/s10107-013-0677-5.

\bibitem{FouGZ} K. Fountoulakis, J. Gondzio and P. Zhlobich,
Matrix-free Interior Point Method for Compressed Sensing Problems, 
Manuscript (2012). arXiv:1208.5435

\bibitem{GonK} C.C. Gonzaga and E.W. Karas, 
Fine tuning Nesterov's steepest descent algorithm
for differentiable convex programming, 
Math. Programming 138 (2013), 141--166.

\bibitem{GonKR} C.C. Gonzaga, E.W. Karas, and D.R. Rossetto,
An optimal algorithm for constrained differentiable convex optimization,
Manuscript (2011).
\url{http://www.optimization-online.org/DB_FILE/2011/06/3053.pdf}

\bibitem{GuLW} M. Gu, L.-H. Lim and C.J. Wu,
PARNES: A rapidly convergent algorithm for accurate recovery of sparse 
and approximately sparse signals,
Numer. Algor. 64 (2013), 321--347.

\bibitem{JudN} A. Juditsky and Y. Nesterov,
Primal-dual subgradient methods for minimizing uniformly convex 
functions,
Manuscript (2010).
\url{http://hal.archives-ouvertes.fr/docs/00/50/89/33/PDF/Strong-hal.pdf}

\bibitem{Lan} G. Lan,
Bundle-level type methods uniformly optimal for smooth and nonsmooth
convex optimization,
Mathematical Programming, (2013) DOI 10.1007/s10107-013-0737-x.

\bibitem{LanLM} G. Lan, Z. Lu and R.D.C. Monteiro,
Primal-dual first-order methods with $O(1/\eps)$ iteration-complexity 
for cone programming,
Math. Programming 126 (2011), 1--29.

\bibitem{MenC} X. Meng and H. Chen,
Accelerating Nesterov's method for strongly convex functions with 
Lipschitz gradient,
Arxiv preprint arXiv:1109.6058 (2011).

\bibitem{NemY} A.S. Nemirovsky and D.B. Yudin,
Problem complexity and method efficiency in optimization,
Wiley, New York 1983.

\bibitem{Nes.dokl} Y. Nesterov,
A method of solving a convex programming problem with convergence 
rate $O(1/k^2)$ (in Russian),
Doklady AN SSSR 269 (1983), 543-547. 
Engl. translation:
Soviet Math. Dokl. 27 (1983), 372--376.

\bibitem{Nes.intL} Y. Nesterov,
Introductory lectures on convex optimization: A basic course,
Kluwer, Dordrecht 2004.

\bibitem{Nes.smooth} Y. Nesterov,
Smooth minimization of non-smooth functions,
Math. Programming 103 (2005), 127--152.

\bibitem{Nes.round} Y. Nesterov,
Rounding of convex sets and efficient gradient methods for linear 
programming problems,
Optim. Methods Softw. 23 (2008), 109--128.

\bibitem{Nes.rel} Y. Nesterov,
Unconstrained convex minimization in relative scale,
Math. Operations Res. 34 (2009), 180--193.

\bibitem{Nes.PD} Y. Nesterov,
Primal-dual subgradient methods for convex problems,
Math. Programming 120 (2009), 221--259.

\bibitem{Ric.imp} P. Richtarik,
Improved algorithms for convex minimization in relative scale,
SIAM J. Optimization 1 (20100), 1141--1167.

\bibitem{Tse} P. Tseng,
On accelerated proximal gradient methods for convex-concave 
optimization,
Manuscript (2008).

\bibitem{YuVGS} J. Yu, S.V.N. Vishvanathan, S. G\"unter and
N.N. Schraudolph,
A quasi-Newton approach to nonsmooth convex optimization problems 
in machine learning,
J. Machine Learning Res. 11 (2010), 1145--1200.



\end{thebibliography}
\end{document}